\documentclass[12pt]{amsart}

\usepackage{verbatim}
\usepackage{amssymb}
\usepackage{upref}
\usepackage[all]{xy}

\usepackage[lite]{amsrefs}

\newtheorem{thm}{Theorem}[section]
\newtheorem*{thm*}{Theorem}
\newtheorem{lem}[thm]{Lemma}
\newtheorem{cor}[thm]{Corollary}
\newtheorem*{cor*}{Corollary}

\theoremstyle{definition}

\newcommand{\thmref}[1]{Theorem~\textup{\ref{#1}}}
\newcommand{\corref}[1]{Corollary~\textup{\ref{#1}}}
\newcommand{\lemref}[1]{Lemma~\textup{\ref{#1}}}

\newcommand{\righttext}[1]{\qquad\text{#1 }}
\renewcommand{\for}{\righttext{for}}

\newcommand{\midtext}[1]{\quad\text{ #1 }\quad}
\renewcommand{\and}{\midtext{and}}
\newcommand{\ie}{\textit{i.e.}}

\newcommand{\new}[1]{{#1}}

\newcommand{\K}{\mathcal K}

\newcommand{\M}{\mathcal M}

\renewcommand{\a}{\alpha}
\renewcommand{\b}{\beta}

\renewcommand{\d}{\delta}
\renewcommand{\l}{\lambda}
\renewcommand{\phi}{\varphi}
\newcommand{\p}{\phi}
\newcommand{\e}{\epsilon}

\newcommand{\s}{\sigma}
\renewcommand{\t}{\theta}
\newcommand{\Chi}{\raisebox{2pt}{\ensuremath{\chi}}}

\renewcommand{\P}{\Phi}
\renewcommand{\L}{\Lambda}

\DeclareMathOperator{\ind}{Ind}
\DeclareMathOperator{\ad}{Ad}

\newcommand{\id}{\text{\textup{id}}}

\renewcommand{\:}{\colon}

\newcommand{\inv}{^{-1}}

\newcommand{\what}{\widehat}

\newcommand{\eg}{\emph{e.g.}}

\begin{document}

\title[Landstad characterization]
{Landstad's characterization for full crossed products}

\author{S. Kaliszewski}
\address{Department of Mathematics and Statistics
\\Arizona State University
\\Tempe, Arizona 85287}
\email{kaliszewski@asu.edu}

\author{John Quigg}
\address{Department of Mathematics and Statistics
\\Arizona State University
\\Tempe, Arizona 85287}
\email{quigg@math.asu.edu}



\date{\today}

\begin{abstract}
Full $C^*$-crossed products
by actions of locally compact groups
are characterized via the
existence of suitable maximal coactions,
in analogy with Landstad's characterization of reduced crossed
products.
\end{abstract}

\maketitle

\section{Introduction} It (almost) goes without saying that two of the
fundamental issues when dealing with $C^*$-dynamical systems are: 
(1)~When can a given $C^*$-algebra be realized as a crossed product 
by a given group? 
(2)~When is a given crossed product isomorphic to a given $C^*$-algebra?
For reduced crossed products by
actions of locally compact groups, 
Landstad (\cite{lan:dual})
answered question~(1) in terms of 
the existence of suitable reduced coactions.
An analogous solution for crossed products by reduced coactions
was given by the second author in \cite{qui:landstad}. 

Answers to question~(2) can be derived from a large class of results
which, in various situations, show that an equivariant homomorphism
between $C^*$-algebras is faithful if it is faithful on a certain
subset of ``fixed points'' in its domain.  Such results exist for
actions of compact groups (folklore), normal coactions of locally 
compact groups
(\cite{qui:discrete}, \cite{qui:landstad}), 
graph $C^*$-algebras (gauge-invariant
uniqueness theorems), and dual coactions on reduced crossed products
(also probably folklore --- see Theorem~\ref{faithful reduced}). 

In this paper we answer questions~(1) and~(2) 
for full crossed products by actions of locally compact groups
(Theorems~\ref{characterize} and~\ref{faithful}, respectively).
In each case, the answer involves the existence of a certain 
maximal coaction; this is natural, since maximal coactions are precisely
those for which full-crossed-product duality holds. 
Also in each case, the answer is gotten by reducing to the situation
of reduced crossed products and normal coactions.  
(See Section~\ref{prelim} for more background on coactions.)  
Although direct proofs almost certainly exist, this 
gives quicker results because the reduced theorems
(\ref{landstad thm} and~\ref{faithful reduced}) are easily
derived from existing work (mostly of Landstad) on reduced coactions.

In Section~\ref{sec:application}, we give an application 
which identifies the crossed product $A^K\times_{\beta}G/K$ of
a fixed-point algebra by a quotient group with a corner of
$A\times_\alpha G$, for $K$ compact.  This result
was motivated by our ongoing study of Hecke $C^*$-algebras
in terms of projections, initiated in \cite{hecke05}. 

\section{Preliminaries}
\label{prelim}

Throughout this paper, $A$ and $B$ denote $C^*$-algebras,
and $G$ denotes a locally compact group.
We adopt the conventions of
\cite{enchilada} for actions and normal coactions, and of \cite{ekq} for maximal coactions.
To avoid ambiguity, we list here a few of our notational conventions.

\subsection*{Actions}
If $\a$ is an action of a locally compact group
$G$ on a $C^*$-algebra $A$, we denote the
reduced crossed product by $A\times_{\a,r} G$, the full
crossed product by $A\times_\a G$,
the canonical covariant homomorphism of $(A,G)$ into the multiplier algebra
$\M(A\times_{\a,r} G)$ by $(i^r_A,i^r_G)$,
the canonical covariant homomorphism of $(A,G)$ into
$\M(A\times_\a G)$ by $(i_A,i_G)$, and the regular representation by
$\L\:A\times_\a G\to A\times_{\a,r} G$.
(Note that $i_A^r = \L_\a\circ i_A$ and $i_G^{r} = \L\circ i_G$.)
When there is more than one action around, we 
will write $i_G^\a$ and $i_G^\b$, for example, 
to distinguish the associated maps.

We identify $G$ with its canonical image in $\M(C^*(G))$, and denote
the left regular representation of $G$ by $\l$ and the compact operators
on $L^2(G)$ by $\K$.

\subsection*{Coactions}
We tacitly assume all our coactions 
are nondegenerate.
If $\d\:A\to \M(A\otimes C^*(G))$ is a coaction of $G$ on $A$,
we denote the crossed product by $A\times_\d G$ and
the canonical covariant homomorphism of $(A,C_0(G))$ into
$\M(A\times_\d G)$ by $(j_A,j_G)$.
The integrated form of the covariant homomorphism
\[
\bigl((\id\otimes\l)\circ\d,1\otimes M)\:
(A,C_0(G))\to \M(A\otimes \K)
\]
(where $M$ is the multiplication representation of $C_0(G)$ on
$L^2(G)$) is faithful.

For every action $(A,\a)$ of $G$ there is a \emph{dual coaction}
$(A\times_\a G,\what\a)$ of $G$, namely the integrated form of the
covariant pair
\begin{align*}
a&\mapsto i_A(a)\otimes 1,\qquad a\in A
\\
s&\mapsto i_G(s)\otimes s,\qquad s\in G.
\end{align*}
There is also a coaction $(A\times_{\a,r} G,\what\a^n)$ of $G$
which is the integrated form of the pair
\begin{align*}
a&\mapsto i^r_A(a)\otimes 1,\qquad a\in A
\\
s&\mapsto i^r_G(s)\otimes \l_s,\qquad s\in G.
\end{align*}
The regular representation
$\L$
is equivariant with respect to $\what\a$ and $\what\a^n$.

For every coaction $(A,\d)$ of $G$ there is a \emph{dual action}
$(A\times_\d G,\what\d)$ of $G$, such that for each $s\in G$ the automorphism
$\what\d_s$ is the integrated form of the covariant pair
\begin{align*}
a&\mapsto j_A(a),\qquad a\in A
\\
f&\mapsto j_G(\rho_s(f)),\qquad f\in C_0(G),
\end{align*}
where $\rho_s(f)(t)=f(ts)$.

The crossed product is functorial; in particular, 
for every equivariant homomorphism
$\p\:(A,\a)\to (B,\b)$
between actions there is a unique equivariant homomorphism
$\p\times G\:(A\times_\a G,\what\a)\to (B\times_\b G,\what\b)$
making the diagram
\[
\xymatrix{
A
\ar[d]_{i_A}
\ar[r]^-\p
&B
\ar[d]^{i_B}
\\
\M(A\times_\a G)
\ar[r]_-{\p\times G}
&\M(B\times_\b G)
}
\]
commute, and for every equivariant homomorphism
$\p\:(A,\d)\to (B,\e)$
between coactions there is a unique equivariant homomorphism 
$\p\times G\:(A\times_\d G,\what\d)\to (B\times_\e G,\what\e)$
making the diagram
\[
\xymatrix{
A
\ar[d]_{j_A}
\ar[r]^-\p
&B
\ar[d]^{j_B}
\\
\M(A\times_\d G)
\ar[r]_-{\p\times G}
&\M(B\times_\e G)
}
\]
commute.

For any coaction $(A,\d)$, the \emph{canonical surjection}
is the map
\[
\P_A
:=\bigl((\id\otimes\l)\circ\d\bigr)\times(1\otimes M)\otimes(1\otimes p)
\]
of $A\times_\d G\times_{\what\d} G$ onto $A\otimes \K$. 
$\d$ is \emph{maximal} if $\P_A$ is faithful, hence an isomorphism.
Thus maximality means that crossed-product duality holds for full
crossed products.

By \cite{ekq}*{Theorem~3.3},  
every coaction $(A,\d)$ has a \emph{maximalization};
\ie, there is a maximal coaction $(A^m,\d^m)$ and an equivariant
surjection
$\vartheta\:(A^m,\d^m)\to (A,\d)$
such that the induced map
$\vartheta\times G\:A^m\times_{\d^m} G\to A\times_\d G$
is an isomorphism.
Moreover, the maximalization is essentially unique in the sense that
if $(C,\e)$ is another maximalization with equivariant surjection
$\varphi$, then there is an equivariant isomorphism $\chi$ making
the diagram
\[
\xymatrix{
(A^m,\d^m)
\ar[r]^-\chi_-\cong
\ar[d]_{\vartheta}
&(C,\e)
\ar[dl]^\varphi
\\
(A,\d)
}
\]
commute.
We remark that the construction of the maximalization in \cite{ekq}
is not functorial --- it involves a choice of a minimal projection in
$\K$. Fischer (\cite{fischer}) has given a functorial version, even for
coactions of quantum groups, but we will not need this here.

A coaction $(A,\d)$ 
\emph{normal} if $j_A$ is faithful; equivalently, if
$(\id\otimes\l)\circ\d$ is faithful.
By \cite{ekq}*{Proposition~2.2}, 
normality is equivalent to crossed-product duality holding for
reduced crossed products; \ie, $\d$ is normal if and only if the
canonical surjection $\P_A$ factors through an isomorphism of the
reduced double crossed product:
\[
\xymatrix{
A\times_\d G\times_{\what\d} G
\ar[r]^-{\P_A}
\ar[d]_\L
&A\otimes \K
\\
A\times_\d G\times_{\what\d,r} G
\ar[ur]_\cong
}
\]
By \cite{qui:fullreduced}*{Proposition~2.6}, every coaction $(A,\d)$
has a \emph{normalization}, \ie, there is a normal
coaction $(A^n,\d^n)$ and an equivariant surjection
$\psi_A\:(A,\d)\to (A^n,\d^n)$
such that the induced map
$\psi_A\times G\:A\times_\d G\to A\times_{\d^n} G$
is an isomorphism.
We have
\[
\ker \psi_A=\ker j_A,
\]
and in fact the construction in \cite{qui:fullreduced} uses
$A^n=j_A(A)$.

Normalization is functorial; in particular, for every
equivariant homomorphism
$\p\:(A,\d)\to (B,\e)$
between coactions there is a unique equivariant homomorphism $\p^n$
making the diagram
\[
\xymatrix{
(A,\d) \ar[r]^-\p \ar[d]_{\psi_A}
&(B,\e) \ar[d]^{\psi_B}
\\
(A^n,\d^n) \ar[r]_-{\p^n}
&(B^n,\e^n)
}
\]
commute, and if $\p$ is an isomorphism then so 
is $\p^n$ (\cite{ekq}*{Lemma~2.1}).

If $(B,\a)$ is an action of $G$, 
then the coaction $(B\times_\a G,\what\a)$ is maximal, 
the coaction $(B\times_{\a,r} G,\what\a^n)$ is normal, 
and the regular representation
\[
\L\:(B\times_\a G,\what\a)\to (B\times_{\a,r} G,\what\a^n)
\]
is both a maximalization of $(B\times_{\a,r} G,\what\a^n)$ 
and a normalization of $(B\times_\a G,\what\a)$ 
(\cite{enchilada}*{Proposition~A.61}, \cite{ekq}*{Proposition~3.4}).


\section{Characterization}
\label{sec:characterize}

We begin with a version, in modern language, of Landstad's
characterization, including a version of his uniqueness clause.

\begin{thm}[Landstad]
\label{landstad thm}
Let $A$ be a $C^*$-algebra 
and let $G$ be a locally compact group. 
There exists an action $(B,G,\alpha)$ 
and an isomorphism $\theta\colon A\to B\times_{\alpha,r}G$
if and only if 
there exists a normal coaction $\delta$ of $G$ on $A$
and a strictly continuous unitary homomorphism 
$u\colon G\to \M(A)$ 
such that
\[
\d(u_s)=u_s\otimes s\quad\text{ for all } s\in G.
\]

Given $\delta$ and $u$ as above, $B$, $\a$, and $\t$ can be chosen
such that $\theta$ is $\delta-\what{\a}^n$ equivariant
and $\theta\circ u = i_G^{\alpha,r}$\textup;  
moreover, if $(C,G,\b)$ is another action 
and $\s\colon A\to C\times_{\beta,r}G$ 
is a $\delta-\what{\b}^n$ equivariant isomorphism
such that $\s\circ u=i^{\beta,r}_G$, 
then there is an isomorphism $\phi\colon (B,\alpha)\to (C,\beta)$ 
such that $(\phi\times_r G)\circ\theta = \sigma$.
\end{thm}

\begin{proof}
We must show that what Landstad actually proved in \cite{lan:dual}
implies the above theorem. First of all, Landstad worked with
\emph{reduced} coactions, which involve the reduced group
$C^*$-algebra $C^*_r(G)$ rather than $C^*(G)$. However, as
shown in \cite{qui:fullreduced}, (nondegenerate) reduced coactions
are in 1-1 correspondence with normal (full) coactions, with the
same crossed products, so we can safely state Landstad's results in
terms of normal coactions. Careful examination of Landstad's proof of
\cite{lan:dual}*{Theorem~3} shows that, given 
\new{$\delta$ and $u$ as in}
our hypotheses,
there exists an action $(B,\a)$ and an equivariant isomorphism
\[
\t\:(A,\d)\xrightarrow{\cong} (B\times_{\a,r} G,\what\a^n)
\]
such that $\t\circ u=i^{\a,r}_G$.
\new{Conversely, given $B$, $\alpha$, and $\theta$
as in our hypotheses, one can easily construct
$\delta$ and $u$ from $\what{\alpha}^n$ and $i_G^{\alpha,r}$ using $\t$.}

However, Landstad's uniqueness clause is stated differently from ours.
Translating into our notation, his uniqueness is in the form
of a characterization of the image of $B$ in the multipliers of the
reduced crossed product.  Specifically, it states that
\begin{align*}
i^r_B(B)
&=\{c\in \M(B\times_{\a,r} G) \mid \what\a^n(c)=c\otimes 1,
\\&\qquad\qquad s\mapsto \ad i^{\a,r}_G(s)(c) \text{ is norm continuous,} 
\\&\qquad\qquad\text{and }ci^{\a,r}_G(f),i^{\a,r}_G(f)c\in B\times_{\a,r} G
\text{ for all $f\in C_c(G)$}\}.
\end{align*}
It follows that
\begin{align*}
\t\inv\circ i^r_B(B)
&=\{a\in \M(A)\mid \d(a)=a\otimes 1,
\\&\qquad\qquad s\mapsto \ad u_s(a) \text{ is norm continuous,}
\\&\qquad\qquad\text{and } au(f),u(f)a\in A
\text{ for all $f\in C_c(G)$}\}.
\end{align*}
We must verify that this implies our uniqueness clause. Suppose we
also have an action $(C,\b)$ and an isomorphism 
$\s\:(A,\d)\xrightarrow{\cong}
(C\times_{\b,r} G,\what\b^n)$ such that $\s\circ u=i^{\b,r}_G$. Then we
have
\begin{align*}
\s\inv\circ i^r_C(C)
&=\{a\in \M(A)\mid \d(a)=a\otimes 1,
\\&\qquad\qquad s\mapsto \ad u_s(a)\text{ is norm continuous,}
\\&\qquad\qquad\text{and } au(f),u(f)a\in A
\text{ for all $f\in C_c(G)$.}\}.
\end{align*}
Thus $\t\inv(i^r_B(B))$ and $\s\inv(i^r_C(C))$ are the same
subalgebra of $\M(A)$, so
\[
\phi:=i^r_C{}\inv\circ \s\circ \t\inv\circ i^r_B
\]
is an isomorphism of $B$ onto $C$; it follows from the definitions
that $\phi$ is equivariant for the actions $\a$ and $\b$\new{, and
that $\phi\times_r G = \sigma\circ\theta^{-1}$}.
\end{proof}

Here is our version of Landstad's characterization for full
crossed products:

\begin{thm}
\label{characterize}
Let $A$ be a $C^*$-algebra 
and let $G$ be a locally compact group. 
There exists an action $(B,G,\alpha)$ 
and an isomorphism $\theta\colon A\to B\times_{\alpha}G$
if and only if 
there exists a maximal coaction $\delta$ of $G$ on $A$
and a strictly continuous unitary homomorphism 
$u\colon G\to \M(A)$ 
such that
\[
\d(u_s)=u_s\otimes s\quad\text{ for all } s\in G.
\]

Given $\delta$ and $u$ as above, 
$B$, $\a$, and $\t$ can be chosen such that
$\theta$ is $\delta-\what{\a}$ equivariant
and $\Lambda_\alpha\circ\theta\circ u = i_G^{\a,r}$,
where $\L_\alpha\:B\times_\a G\to B\times_{\a,r} G$ 
is the regular representation\textup;
moreover, 
if $(C,G,\b)$ is another action 
and $\s\colon A\to C\times_{\beta}G$ 
is a $\delta-\what{\b}$ equivariant isomorphism
such that $\Lambda_\beta\circ\s\circ u=i^{\b,r}_G$, 
then there is an isomorphism $\phi\colon (B,\alpha)\to (C,\beta)$ 
such that 
$\Lambda_\beta\circ(\phi\times G)\circ\theta = \Lambda_\beta\circ\sigma$.
\end{thm}

\begin{proof}
\new{For the forward implication, $\delta$ and $u$ are easily
constructed from $\what\a$ and $i_G^\a$ using $\t$.  
For the reverse implication,}
let $\psi\:(A,\d)\to (A^n,\d^n)$ be the normalization.
Then $u^n:=\psi\circ u\:G\to \M(A^n)$ 
is a strictly continuous unitary homomorphism such that
\[
\d^n(u^n_s)= u^n_s\otimes s,
\]
so by \thmref{landstad thm} 
there exists an action $(B,\a)$ 
and an isomorphism we will call
$\theta^n\:(A^n,\d^n)\to (B\times_{\a,r} G,\what\a^n)$ 
such that $\theta^n\circ u^n=i^{\a,r}_G$.
Note that 
\[
\theta^n\circ \psi\:(A,\d)\to (B\times_{\a,r} G,\what\a^n)
\]
is maximalization of $(B\times_{\a,r} G,\what\a^n)$. 

On the other hand,
$\L_\alpha\:(B\times_\a G,\what\a)\to (B\times_{\a,r} G,\what\a^n)$
is also a maximalization of $(B\times_{\a,r} G,\what\a^n)$. 
Thus, by uniqueness of maximalizations, 
there exists an isomorphism $\t$ making the diagram
\[
\xymatrix{
(A,\d) \ar[r]^-\t \ar[d]_\psi
&(B\times_\a G,\what\a) \ar[d]^{\L_\alpha}
\\
(A^n,\d^n) \ar[r]_-{\theta^n}
&(B\times_{\a,r} G,\what\a^n)
}
\]
commute, and we have
$\L_\alpha\circ \t\circ u
=\theta^n\circ \psi\circ u
=\theta^n\circ u^n =i^{\a,r}_G$.

For the uniqueness, 
suppose we have another action $(C,G,\b)$ 
and an isomorphism $\s\:(A,\d)\xrightarrow{\cong} (C\times_\b G,\what\b)$
such that $\L_\beta\circ\s\circ u=i^{\b,r}_G$.
Since the regular representation
$\L_\beta\:(C\times_\b G,\what\b)\to (C\times_{\b,r} G,\what\b^n)$
is a normalization, 
by functoriality of normalization 
we have an isomorphism $\s^n$ making the
diagram
\[
\xymatrix{
(A,\d) 
\ar[r]^-\s_-\cong
\ar[d]_\psi
&(C\times_\b G,\what\b) 
\ar[d]^{\L_\beta}
\\
(A^n,\d^n) 
\ar[r]_-{\s^n}^-\cong
&(C\times_{\b,r} G,\what\b^n)
}
\]
commute. Thus (with $u^n=\psi\circ u$ as above) we have
\[
\s^n\circ u^n
= \s^n\circ\psi\circ u
= \L_\beta\circ\s\circ u
= i_G^{\beta,r},
\]
so by \thmref{landstad thm}
there is an isomorphism $\phi\colon(B,\a)\xrightarrow{\cong}(C,\b)$
\new{such that 
\[
\L_\b\circ(\phi\times G)\circ\t
= (\phi\times_r G)\circ\t^n
= \sigma^n = \Lambda_\beta\circ\sigma.
\]
}

\end{proof}

\section{Fidelity}
\label{sec:faithful}

We begin with a fidelity criterion for homomorphisms of reduced
crossed products. This result is probably folklore; we could not find
it in the literature, so we include a proof for completeness.

\begin{thm}
\label{faithful reduced}
Let $(B,\a)$ be an action of a locally compact group $G$,
and let $\phi\colon B\times_{\a,r}G\to A$ be a surjective
homomorphism.
Then $\phi$ is an isomorphism 
if and only if $\p\circ i^r_B$ is faithful 
and there exists 
a normal coaction $\delta$ of $G$ on $A$ such that 
$\phi$ is $\what{\a}^n-\delta$ equivariant.
\end{thm}

\begin{proof}

\new{For the forward implication, $\phi\circ i_B^r$ is clearly faithful,
and $\delta$ is easily constructed from $\what{\a}^n$ using $\phi$.}
We adapt to our purposes a trick of Landstad \cite{lan:dualcompact}
\new{for the reverse implication}.
Let 
$(\pi,u)$ be a covariant homomorphism of $(B,G,\alpha)$
into $\M(A)$ such that 
\[
\pi\times u = \p\circ\L_\a\:B\times_a G\to A.
\]
Since the coaction $\d$ is normal, the homomorphism
\[
\tau:=(\id\otimes\l)\circ\d\:A\to \M(A\otimes \K)
\]
is faithful. 
Elementary computations using equivariance of $\pi\times u$ with respect to $\what\a$ and $\d$ show that
\[
\tau\circ (\pi\times u)
=((\pi\times u)\otimes \l)\circ\what\alpha
= (\pi\otimes 1)\times(u\otimes\l),
\]
which is unitarily equivalent to the homomorphism $\ind\pi$ induced from
$\pi$ (see, \eg, \cite{enchilada}*{Lemma~A.18}). 
Since 
\[
\pi=(\pi\times u)\circ i_B = \phi\circ\Lambda_\alpha\circ i_B
= \phi\circ i_B^r
\]
is faithful by hypothesis,
$\ind\pi$ is a regular representation. Thus
\[
\tau\circ(\pi\times u)=\tau\circ\p\circ\L_\a
\]
has the same kernel as $\L_\a$, so $\tau\circ\p$ is faithful, hence so is $\p$.

\end{proof}

In order to deduce an analogous fidelity result for full crossed products,
we will need the following lemmas:

\begin{lem}
\label{fixed}
Let $(A,G,\d)$ be a coaction, with normalization
\[
\xymatrix{
(A,\d) \ar[r]^-{\psi}
&(A^n,\d^n).
}
\]
Then $\psi$ is 1-1 on
\[
\M(A)^\d:=\{a\in \M(A)\mid \d(a)=a\otimes 1\}.
\]
\end{lem}

\begin{proof}
By \cite{enchilada}*{Proposition~A.61}, we can take
$\psi=(\id\otimes\l)\circ\d$,
so that for $a\in \M(A)^\d$ we have
\[
\psi(a)=(\id\otimes\l)(\d(a))
=(\id\otimes\l)(a\otimes 1)
=a\otimes 1,
\]
which is $0$ if and only if $a$ is.
\end{proof}

\begin{lem}
\label{isomorphism}
Let $(C,G,\e)$ and $(A,G,\d)$ be maximal coactions, and let
\[
\xymatrix{
(C,\e) \ar[r]^-\p & (A,\d)
}
\]
be a surjective equivariant homomorphism. If the normalization $\p^n$
is an isomorphism, then so is $\p$.
\end{lem}

\begin{proof}
By \cite{ekq}*{Lemma~2.1}, we have a commutative diagram
of normalizations
\[
\xymatrix{
C \ar[r]^-\p \ar[d]_{\psi_C}
&A \ar[d]^{\psi_A}
\\
C^n \ar[r]_-{\p^n}
&A^n.
}
\]
Applying functoriality of crossed products,
we get a commutative diagram
\[
\xymatrix
@C+50pt
@R+10pt
{
C\otimes\K
\ar[r]^-{\p\otimes\id}
&A\otimes\K
\\
C\times_\e G\times_{\what\e} G
\ar[u]^{\P_C}_\cong
\ar[r]^-{\p\times G\times G}
\ar[d]_{\psi_C\times G\times G}^\cong
&A\times_\d G\times_{\what\d} G
\ar[u]_{\P_A}^\cong
\ar[d]^{\psi_A\times G\times G}_\cong
\\
C^n\times_{\e^n} G\times_{\what{\e^n}} G
\ar[r]^-{\p^n\times G\times G}
&A^n\times_{\d^n} G\times_{\what{\d^n}} G,
}
\]
where the lower two vertical arrows are isomorphisms by definition of
normalization, and the upper two vertical arrows (the canonical
surjections) are isomorphisms because $\e$ and $\d$ are maximal. 
Assuming that $\p^n$ (and hence $\p^n\times G\times G)$ is an isomorphism, 
it follows that $\p\otimes\id$ is an isomorphism, and therefore so is $\p$.
\end{proof}

\begin{cor}
\label{faithful}
Let $(B,\a)$ be an action of a locally compact group $G$,
and let 
$\p\:B\times_\a G\to A$ be a surjective
homomorphism.
Then $\phi$ is an isomorphism 
if and only if $\p\circ i_B$ is faithful 
and there exists a maximal coaction $\delta$ of $G$ on $A$
such that $\phi$ is $\what\a-\delta$ equivariant.
\end{cor}

\begin{proof}
\new{The forward implication is straightforward. For the reverse
implication,}
by \cite{ekq}*{Lemma~2.1}, 
the equivariant homomorphism $\p$ fits into a commutative diagram
\[
\xymatrix{
(B\times_\a G,\what\a) \ar[r]^-\p \ar[d]_{\L_\a}
&(A,\d) \ar[d]_{\psi}
\\
(B\times_{\a,r} G,\what\a^n) \ar[r]_-{\p^n}
&(A^n,\d^n),
}
\]
where the vertical arrows are normalizations. 
It is easily verified that
$\p\circ i_B(B) \subseteq \M(A)^\delta$,
so $\psi$ is 1-1 on the range of $\p\circ i_B$
by \lemref{fixed}.
Since $\p\circ i_B$ is faithful, so is 
\[
\p^n\circ i^r_B = \p^n\circ\L_\a\circ i_B = \psi\circ\p\circ i_B. 
\]
Thus \thmref{faithful reduced} tells us
that $\p^n$ is an isomorphism, and it follows from 
\lemref{isomorphism} that $\p$ is also an isomorphism.
\end{proof}


\section{Application}
\label{sec:application}

We apply \corref{faithful}
to the following situation: let $\a$ be
an action of $G$ on $A$, and let $K$ be a compact normal subgroup
of $G$. Normalize the Haar measure on $G$ so that $K$ has measure
$1$, making the characteristic function $\Chi_K$ a projection
in $C^*(G)$.
Let $A^K$ denote the fixed-point subalgebra of $A$
under $K$. Then the action of $G$ leaves $A^K$ invariant, so
there is a unique action $\b$ of $G/K$ on $A^K$ such that
$\b_{sK}(a)=\a_s(a)$ for all $a\in A^K$. 
Let $p=i_G(\Chi_K)$, which is a projection in $\M(A\times_\a G)$
which commutes with $i_A(A^K)$.
(Recall that $(i_A,i_G)$ denotes
the canonical covariant homomorphism of $(A,G,\a)$ into
$\M(A\times_\a G)$.)
Routine calculations show that the maps
\begin{align*}
\s\:A^K&\to A\times_\a G&\tau\:G/K&\to \M(A\times_\a G)\\
\s(a)&=i_A(a)p&\tau(sK)&=i_G(s)p
\end{align*}
give a covariant homomorphism $(\s,\tau)$ of $(A^K,G/K,\b)$ to
$\M(A\times_\a G)$.
Moreover, the integrated form $\t:=\s\times \tau$ is a surjection of
$A^K\times_\b G/K$ onto $p(A\times_\a G)p$.

\begin{cor}
With the above hypotheses and notation, the map
\[
\t\:A^K\times_\b G/K\to p(A\times_\a G)p
\]
is an isomorphism.
\end{cor}

\begin{proof}
It remains to show that $\t$ is injective, and to do this
we aim to apply \corref{faithful}.
The dual coaction $\what\a$ of $G$ on $A\times_\a G$
is maximal by \cite{ekq}*{Proposition~3.4}.
Moreover, $\what\a$
restricts to a coaction
$\what\a|:=(\id\otimes q)\circ\what\a$ of $G/K$ on $A\times_\a G$
(where $q\:C^*(G)\to C^*(G/K)$ denotes the quotient homomorphism),
and
by \cite{kq:fullmans}*{Corollary~7.2} this coaction is also maximal.

The projection $p$ is invariant under the coaction $\what\a$,
hence under $\what\a|$, so by the elementary lemma
below the
restriction $\d$ of $\what\a|$ to the corner $p(A\times_\a G)p$
is a maximal coaction.
Routine calculations verify that
$\d\circ\t=(\t\otimes\id)\circ\what\b$,
so $\t$ is $\what\beta-\delta$ equivariant. 
$\t\circ i_B=\sigma$ is evidently injective,
so $\theta$ is injective by \corref{faithful}.
\end{proof}

\begin{lem}
\label{projection}
Let $\d$ be a maximal coaction of $G$ on $A$,
and let $p$ be an invariant projection in $\M(A)$.
Then $\e\:=\d|_{pAp}$ is a maximal coaction of $G$ on
the corner $pAp$.
\end{lem}

\begin{proof}
Let $B:=pAp$.  
Routine calculations show that $\e$ is a coaction of $G$ on $B$.
For maximality, 
as in the proof of \cite{ekq}*{Proposition~3.5}, 
the canonical surjection
$\Phi_B\colon B\times_\e G\times_{\what{\e}} G
\to B\otimes \K$ can be identified with a restriction
of $\Phi_A$, and therefore $\Phi_B$ is injective
because $\Phi_A$ is.
\end{proof}

\end{document}